\input amstex
 \documentstyle{amsppt}
 \magnification=\magstep1
 \vsize=24.2true cm
 \hsize=15.3true cm
 \nopagenumbers\topskip=1truecm
 \headline={\tenrm\hfil\folio\hfil}

 \TagsOnRight

\hyphenation{auto-mor-phism auto-mor-phisms co-homo-log-i-cal co-homo-logy
co-homo-logous dual-izing pre-dual-izing geo-metric geo-metries geo-metry
half-space homeo-mor-phic homeo-mor-phism homeo-mor-phisms homo-log-i-cal
homo-logy homo-logous homo-mor-phism homo-mor-phisms hyper-plane hyper-planes
hyper-sur-face hyper-sur-faces idem-potent iso-mor-phism iso-mor-phisms
multi-plic-a-tion nil-potent poly-nomial priori rami-fication sin-gu-lar-ities
sub-vari-eties sub-vari-ety trans-form-a-tion trans-form-a-tions Castel-nuovo
Enri-ques Lo-ba-chev-sky Theo-rem Za-ni-chelli in-vo-lu-tion Na-ra-sim-han Bohr-Som-mer-feld}

 % reduced subscheme
 % minimal surface
\define\rest#1{_{\textstyle{|}#1}} % restriction of map to subset

\define\Span#1{\left<#1\right>} % <x> span or hull of x

 % Volume
 % roman differential dx, dy
 % identity map

% Bourbaki letters
\define\C{\Bbb C} % complex numbers
\define\R{\Bbb R} % real numbers
\define\Z{\Bbb Z} % integers

\define\proj{\Bbb P} % projective space
 % projective plane
 % dual projective plane

% Script letters
 % sheaf of algebras A
\define\sB{{\Cal B}} % conductor ideal
 % sheaf G
 % sheaf J
 % moduli space Mg
\define\Oh{{\Cal O}} % structure sheaf

% short Greeks
\define\al{\alpha}
\define\be{\beta}
\define\de{\delta}

\define\ga{\gamma}

\define\om{\omega}
\define\si{\sigma}
\define\De{\Delta}
\define\Ga{\Gamma}

\define\Om{\Omega}
\define\la{\lambda}

% \mathops
 % adjoint bundle

 % cardinality
 % codimension
 % cokernel of a map

 % image of a map
 % rank of a map

 % modulo
 % base locus of linear system

 % Grassmann variety
 % Hilbert scheme
 % Hom group
 % End group

 % Picard scheme
 % singular locus
 % Stabiliser group
 %second exterior power of a module

% Lie groups

 % unitary group

 \document

  \topmatter
  \title ALAG - quantization
            \endtitle
  \author Nik. Tyurin                     \endauthor

   \address MIIT
   \endaddress
  \email
    tyurin\@tyurin.mccme.ru  ntyurin\@newton.kias.re.kr
   \endemail

\abstract This paper is a continuation of [10]. Here we present
the application of ALAG - programme, introduced in [2], [9], to
the  geometric quantization. The proposed approach is following to
"geometrical formulation of quantum mechanics" ([1]). For our ALAG
- quantization the correspondence principle is satisfied.
\endabstract

\endtopmatter

\head \S 0. Introduction   \endhead

The framework of geometric quantization is usually described as
follows ( see f.e. [5]). Let $(M, \om)$ is a symplectic manifold
which represents a classical mechanical system with finite number
of degrees of freedom. Then one understands the geometric
quantization as a procedure relating with  $(M, \om)$ a Hilbert
space $\Cal H$ together with a correspondence
$$
Q: \C^{\infty}(M \to \R) \to O(\Cal H)
$$
where $O(\Cal H)$ is the space of self adjoint operators which
should satisfy the correspondence principle
$$
Q(\{ f, g \}_{\om}) = \imath \hbar [Q(f), Q(g)]
$$
as well as the irreducibility condition. Unfortunately such
construction couldn't exist (by the van Hove theorem).
Nevertheless one usually uses two basic variants of this
procedure: Souriau - Kostant quantization which  doesn't satisfy
the irreducibility condition and Berezin quantization which
doesn't satisfy the correspondence principle.

But now the question of quantization can be refined  in  more
invariant fashion: we are looking for a quantum mechanical system
which is  closely related to the original classical mechanical
system. Turning to the language of "geometrical formulation of
quantum mechanics" (see [1]) one could reformulate the subject of
geometric quantization as follows. For any symplectic manifold
$(M, \om)$ one would like to construct a K\"ahler manifold $\Cal
K$ (finite or infinite dimensional) endowed with a K\"ahler triple
$(G, I, \Om)$ where $\Om$ is a symplectic form, $I$ is an
integrable complex structure and $G$ is the corresponding
riemannian metric together with a correspondence
$$
\tilde{Q}: C^{\infty}(M \to \R) \to \Cal N \subset C^{\infty}(\Cal
K \to \R)
$$
where $\Cal N$ consists of all functions whose hamiltonian vector
fields preserve riemannian metric $G$:
$$
\Cal N = \{ F \in C^{\infty}(\Cal K \to \R) \quad | \quad
Lie_{X_F} G = 0.\}
$$
One requires for this construction that
$$
\tilde{Q}(\{f, g\}_{\om}) = \{ \tilde{Q}(f), \tilde{Q}(g) \}_{\Om}
\tag *
$$
which means that the correspondence principle is satisfied (see
[1] for explanation why it's equivalent to the last one).

Now one can mention that the so called ALAG - programme (see [2],
[3], [9]) endows us with a candidate --- for any compact smooth
orientable symplectic $2n$ - dimensional manifold $M$ with
integer symplectic form $\om$ one can almost canonically
construct a set of infinite dimensional K\"ahler manifolds
$\sB^{hw,t}_S$ labeled by "integer"parameter $S \in H_n (M, \Z)$
and continuous parameter $t \in \R$. These are called moduli
spaces of half weighted Bohr --- Sommerfeld lagrangian cycles of
fixed volume and topological type.

On the other hand it was pointed out ([9], [10]) that for any such
moduli space $\sB^{hw,t}_S$ one has a natural correspondence: for
every smooth function $f$ on $M$ there is a corresponding
function $F_f$ on $\sB^{hw,1}_S$; these induced functions we will
call special functions. The main goal of the present paper is to
prove the following statement.

\proclaim{Main Theorem} In the setup of ALAG

1) the correspondence
$$
f \mapsto F_f
$$
is an inclusion;

2) for every $f \in C^{\infty}(M \to \R)$ the corresponding
special function $F_f$ satisfies the condition mentioned in the
framework of the geometrical formulation
$$
Lie_{X_{F_f}} G = 0;
$$

and

3) the correspondence principle in form (*) is satisfied.
\endproclaim

Thus one gets a new approach to the geometric quantization
construction via ALAG - programme. At the same time this approach
gives us the first application of ALAG - programme in pure
geometrical setup. But really the spectrum of  validity of ALAG
is much more wider and one hopes to see a number of new results
in a nearest future.

I would like to cordially thank Korean Institute of Advanced Study
for hospitality during this work. This work would be impossible
without uninterrupted dialog with the authors of ALAG took place
all the time we were staying in Seoul.

\head \S 1. Geometrical formulation of quantum mechanics
\endhead

In this subsection we follow [1] preserving the notations.

Let us briefly recall the main ideas of the  geometrical
formulation.

In classical mechanical setup one has a symplectic manifold $M$
which represents the corresponding phase space such that states
are represented just by points of $M$. The set of classical
observables is the space of smooth functions $C^{\infty} (M \to
\R)$ endowed with a Poisson algebra structure which is given by
the symplectic form. The measurement aspects are quite simple: one
just takes the volume of observables over a state $x \in M$.
Every smooth function $f$ generates the corresponding hamiltonian
vector field
$$
X_f = \om^{-1}(df)
$$
which preserves the symplectic form
$$
Lie_{X_f} \om = 0
$$
and induces a motion over $M$ which can be understood as a
continuous family of  symplectomorphisms  of $M$ and its  germ
corresponds to $X_f$.

In quantum mechanics usually one deals with a Hilbert space $\Cal
H$ so a complex vector space (finite or infinite dimensional)
endowed with a hermitian inner product $<.,.>$. Quantum states are
represented by rays in $\Cal H$. Quantum observables are self
adjoint operators such that for every $A, B$ one has a skew
symmetric operation
$$
\hat A, \hat B \mapsto \imath \hbar [\hat A, \hat B]
$$
which endows the space of self adjoint operators with a Lie
algebra structure. The dynamics of the system is described by
Schr\"odinger equation
$$
\psi = - \frac{\imath}{\hbar} \hat H \psi
$$
where $\hat H$ is a preferred observable called the hamiltonian of
the system. The corresponding flow is generated by 1 - parameter
group $exp (\imath \hat H t)$. The measurement aspects are much
more complicated: the ideal measurement of an observable $\hat A$
in a state $\psi \in \Cal H$ yields an eigenvalue of $\hat A$ and
immediately after the measurement the state is thrown into the
corresponding eigenstate. So roughly speaking for $\hat A$ and
$\psi$ there is a set of  projections of the state $\psi$ to the
eigenstates of $\hat A$ which one understands as probabilities
(amplitudes)(one can twist every orthonormal eigenbasis such that
all coordinates are real nonnegative numbers).

The starting point of discussion is the fact that the space of
physical states is not the Hilbert space $\Cal H$ itself but the
space of rays in $\Cal H$. So instead of $\Cal H$ one can consider
its projectivization
$$
\Cal P = \proj (\Cal H)
$$
which is a K\"ahler manifold (infinite or finite dimensional as
well as $\Cal H$) endowed with

1) symplectic form $\Om$

2) integrable complex structure $I$

3) the corresponding riemannian metric $G$.

One can translate completely all notions of QM to the language of
this projective space without any references to the original
Hilbert space.

Here it is the vocabulary:

\subheading{Quantum states} Points of $\Cal P$ represent quantum
states.

\subheading{Quantum observables} Instead of self adjoint
operators over $\Cal H$ one takes the corresponding expectation
value functions:
$$
\hat A \mapsto A(\psi) = \frac{1}{2 \hbar} G (\psi, \hat A \psi).
$$
Every such function admits well defined restriction to unit sphere
in $\Cal H$ and further to $\Cal P$.

The complement question is: how one can recognize is a smooth
function $F$ over $\Cal P$ induced by some self adjoint operator
over $\Cal H$? The answer is quite definite: $F$ is induced by a
self adjoint operator  if and only if the hamiltonian vector
$X_F$ preserves the riemannian metric $G$:
$$
Lie_{X_F} G = 0.
$$
Thus one can distinguish  such functions over $\Cal P$ and collect
these in the set $\Cal N$ of special functions.

\subheading{Schr\"odinger vector fields} For  such a function $A:
\Cal P \to \R$ one gets that the hamiltonian vector field $X_A$
on $\Cal P$ is precisely the same as the corresponding
Shroedinger vector field restricted to unit sphere and factorized
by the phase rotations.

\subheading{Commutator} For two observables $\hat A, \hat B$ one
takes the commutator $\hat C = \frac{1}{\imath \hbar} [\hat A,
\hat B]$. Then the induced expectation volume function just
equals to the Poisson bracket
$$
C = \{ A, B \}_{\Om}.
$$

\subheading{Schr\"odinger equation} For a preferred hamiltonian
$\hat H$ over $\Cal H$ the induced motion on $\Cal P$ is precisely
the same as that one which induced by hamiltonian $H$. It means
that Schr\"odinger equation reads as Hamilton equation over $\Cal
P$.

\subheading{Uncertainties} For two special functions $A, B \in
\Cal N$ one has riemannian bracket:
$$
(A, B)_G = \frac{\hbar}{2} G(X_A, X_B).
$$
Then the squared uncertainty of $A$ at the quantum state $p$
equals to
$$
(\De A)^2 (p) = (A, A)_G(p).
$$
The uncertainty relation for two quantum observables looks like
$$
(\De A)(\De B) \geq (\frac{\hbar}{2} \{A, B \}_{\Om})^2 + (A,
B)_G^2.
$$
If in a point $p$ two hamiltonian vector fields $X_A, X_B$ are
$G$ - orthogonal one gets the standard uncertainty relation
$$
(\De A)(\De B) \geq (\frac{\hbar}{2} \{ A, B \}_{\Om}).
$$

\subheading{Eigenstates and eigenvalues}  For a self adjoint
operator $\hat A$ with eigenstates $\psi_i$ and eigenvalues $\la_i
\in \R$ one gets that the corresponding special function $A$ has
these $p_i = \proj (\psi_i)$ as critical points with critical
values
$$
A(p_i) = \la_i.
$$

\subheading{Measurement process. Discrete spectrum} Let one take a
state $\psi$ with respect to orthonormal basis $\psi_i$ consists
of eigenstates of $\hat A$. Then the corresponding amplitudes
$\al(\psi, \psi_i)$ can be found over $\Cal P$ as follows. For $p
= \proj (\psi)$ and $p_i = \proj (\psi_i)$ consider geodesic
distances
$$
d(p, p_i)
$$
with respect to riemannian metric $G$. Then one has the following
equality
$$
\al(\psi, \psi_i) = cos^2 (\frac{d(p, p_i)}{\sqrt{2 \hbar}}).
$$
We already know what are eigenstates and eigenvalues.

\subheading{Measurement process. Continuous spectrum} Almost the
same story just instead of points one takes the corresponding
critical subsets and measures geodesic distance between point $p$
and such a subset.

Thus one could understand the main difference between classical
mechanics and quantum mechanics in the following style (see [1]):
together with a symplectic structure one has a compatible
riemannian structure on the phase space in the last case.

\subheading{Remark} Here let us note that the projectivization
$\proj \Cal H$ is a K\"ahler manifold of special type namely it
is an algebraic manifold. This means that one has over $\proj
\Cal H$ the positive line bundle $\Oh (1)$ which first Chern
class is represented by the K\"ahler form $\om$. The space of
holomorphic sections $H^0(\Oh(1))$ is isomorphic to $\Cal H$ and
in the next section we discuss how one can reconstruct all data
on this Hilbert space.

\head \S 2. Examples
\endhead

Now let us consider two well known examples of geometric
quantization via this projectivization approach.

First of all turn to Souriau - Kostant method. Here we have a
prequantization quadruple $(M, \om, L, a)$ where $L$ is a
hermitian line bundle defined topologically by the first Chern
class $$ c_1(L) = [\om]$$ and $a$ is a hermitian connection
defined by $$ F_a = 2 \pi \imath \om. $$ Then the Hilbert space is
 $$\Cal H = L^2(M, L)$$
--- the completion of the space of all smooth sections of $L \to
M$ with respect to the natural norm
$$
\Vert \si \Vert^2 = \int_M <\si, \si> d \mu,
$$
where $d \mu$ is the Liouville volume form. Then for every
function $f \in C^{\infty}(M \to \R)$ one takes the operator
$$
Q_f: \Cal H \to \Cal H, \quad \quad Q_f = \nabla_{X_f} + 2 \pi
\imath f
$$
acting on $\Cal H$. With  well known relationships in differential
geometry in hands one has
$$
[Q_f, Q_g] = Q_{\{f, g \}_{\om}}.
$$
One can explain this correspondence using some background ideas
which could be called as "dynamical correspondence". The point is
that one has over $M, \om$ the space of symplectomorphism germs
represented by global hamiltonian vector fields over $M$. If one
would like to extend a germ generated by $X_f$ over $M$ to a one
parameter family of automorphisms of the hermitian line bundle $L$
with fixed hermitian connection $a$ then there is unique lifting
of such type which is described exactly by the vector field
corresponds to $Q_f$. More precisely every hamiltonian vector
field $X_f$ induces a germ of linear transformation of our $\Cal
H$ and this transformation is just the exponential of $Q_f$. So
in style of [1] we will write a vector field on $\Cal H$ which
reads just as
$$
Y(\psi) = Q_f \psi.
$$
If one takes "expectation value" of $Q_f$ it would be a pure
imaginary function since  $Q_f$ is a skew hermitian one. One can
rescale this function getting a smooth real function
$\tilde{Q}_f$ on the projective space $\Cal P = \proj (\Cal H)$.

Now every hamiltonian vector field $X_f$ generates a smooth
vector field over $\Cal P$ which preserves both the symplectic
and the riemannian structures over $\Cal P$. This vector field we
will denote as $\Theta(f)$, emphasizing that it is generated by
the smooth "source" function $f$ on the "source" manifold $M$.
This map
$$
\Theta: C^{\infty}(M \to \R) \to Vect( \Cal P)
$$
we will call the {\it dynamical correspondence}.  It plays the
fundamental role in quantization procedure. Obviously the
dynamical correspondence  map satisfies the correspondence
principle
$$
\Theta (\{f, g \}_{\om}) = [\Theta(f), \Theta(g)]
$$
where at the right hand side one has the commutator of two vector
fields. It reflects just the fact that $\Theta$ is defined
directly from dynamical properties of the system.

Generalizing this idea let us say that if one has any object over
$M, \om$ preserved by all symplectomorphisms of $M$ then there
exists the same "dynamical correspondence" between hamiltonian
vector fields and infinitesimal automorphisms of this object. We
will see further how this idea works in ALAG.

Coming back to Souriau - Kostant quantization one gets the
following fact: the induced "dynamical correspondent" vector field
$\Theta(f)$ over $\Cal P$ is equal to Killing vector field of
smooth real function $\tilde{Q}_f$. This claim is based on the
fact that $Q_f$ is a  skew hermitian operator over $\Cal H$. From
the algebraic point of view one can correct the correspondence
$$
f \mapsto Q_f
$$
such that
$$
f \mapsto - \imath Q_f
$$
(we omit the Planck constant every time to clarify the
mathematical aspects of the problem), getting self adjoint
operators and keeping the correspondence principle but the
difficulties which one meets further in this way should be
explained with respect to this projection to "dynamical
correspondence" background. Therefore geometrically the well
known reducibility of the Souriau - Kostant representation has a
root in this "dynamical discrepancy" when $\Theta(Q)$ doesn't
coincide with the hamiltonian vector field $X_{\tilde{Q}_f}$.

Turning to Berezin quantization one should say that it was
originally described in the language of projectivization. Here one
needs an additional structure - an appropriate complex
polarization which can be understood as an integrable complex
structure $I$ compatible with $\om$. Thus we come to the
framework of the algebraic geometry. Then one has a Hilbert
subspace in $\Cal H$ consists of holomorphic sections of $L$
$$
H^0(M_I, L) \subset \Cal H.
$$
Let us suppose that our holomorphic line bundle $L$ is very ample.
Then there is smooth holomorphic inclusion
$$
\psi: M_I \to \proj H^0(M_I, L)^*
$$
defined as usual by the  corresponding complete linear system.
Combining with hermitian conjugation one gets a smooth
antiholomorphic map
$$
\psi_0: M_I \to \proj H^0 (M_I, L),
$$
and its image is called the space of coherent states. Originally
in the Berezin works it was proposed that one takes the space of
symbols  over the image $\psi_0(M)$ which are nothing but the
expectation functions restricted to $\psi_0(M)$. The extensions
of these symbols over all $\proj H^0(M_I, L)$ can be derived
analytically. But we would like to perform this story using so
called Berezin - T\"oplitz operators over $H^0(M_I, L)$. Namely
since our ambient space $\Cal H$ is a Hilbert space then one can
consider Sz\"oge projector
$$
S: \Cal H \to H^0(M_I, L)
$$
which is just the orthogonal projector. Then for every smooth
function $f $ one has the following combination
$$
T_f: H^0(M_I, L) \to H^0(M_I, L), \quad \si \mapsto S(f \si).
$$
These operators are self adjoint. Turning to the geometrical
formulation one gets a  smooth real function $\tilde{T}_f$ over
$\proj H^0(M_I, L)$ which is a function of special type described
in [1]: its hamiltonian vector field should preserve the
riemannian metric (and the complex structure) over $\proj
H^0(M_I, L)$ as well as the symplectic structure. One can compute
this function using Fourier - Berezin  transform. Namely one has a
universal function $u(x, s)$ on the direct product
$$
M \times \proj H^0(M_I, L).
$$
If we represent $\proj H^0(M_I, L)$ by holomorphic sections of
unit norms
$$
\int_M <s, s> d \mu = 1
$$
then
$$
u(x, s) = <s(x), s(x)>
$$
is a nonnegative smooth function. Then the transform reads as
$$
\tilde{T}_f = \int_M f(x) u(x, s) d \mu. \tag 1
$$
On the other hand one couldn't define an analogy of the dynamical
correspondence here. The point is that a generic function $f \in
C^{\infty}(M \to \R)$ generates flow which doesn't preserve the
space $H^0(M_I, L)$ and consequently the projective space $\proj
H^0(M_I, L)$. Thus in the case one should restrict the
consideration to the subspace $\Cal M \subset C^{\infty}(M \to
\R)$ of functions whose flows preserve the projective space $\proj
H^0(M_I, L)$.  For this case we have a brief and precise
description of Berezin quantization contained in [6]. In this
description one combines Souriau - Kostant quantization with
Berezin approach (see [6]). Then it's reasonable to define an
analogy of the dynamical correspondence which is a map
$$
\Theta': \Cal M \to Vect (\proj H^0(M_I, L)).
$$
If $f$ belongs to $\Cal M$ then it induces a special function
$\tilde{Q}_f$ (just the restriction of the "Souriau - Kostant" -
function) such that its Killing vector field coincides with the
induced by the dynamical property
$$
\Theta'(f) = K(\tilde{Q}_f).
$$
And again the dynamical correspondence works establishing that
correspondence principle is satisfied in this case in terms of
"quantum" Poisson bracket. But this approach is rather special
since it allows to quantize only functions of very special type.
In Berezin - T\"oplitz approach one has much more wider situation
but of course the integral operator (1) described above has  big
kernel. At the same time it's much more complicated problem to
define in this approach any kind of the dynamical correspondence.
At least we have not now at hands any simple description of that
one. Because of this we can not prolong the story just mentioning
that in Berezin - T\"oplitz quantization the correspondence
principle holds only asymptotically (see [2]).

\head \S 3. Dynamical correspondence in ALAG
\endhead

Let $(M, \om)$ is a compact smooth orientable $2 n$ - dimensional
symplectic manifold with integer symplectic form $\om$. Then
there exists an infinite dimensional K\"ahler manifold
$\sB^{hw,1}_S$ where $S \in H_n(M, \Z)$ is a homological class
over $M$ and we consider the half weighted cycles of volume 1
just for simplicity (all definitions and constructions are
contained in [2], [3] and in [9], where they were firstly
introduced). Moreover this moduli space is an algebraic manifold
since the K\"ahler form represents the first Chern class of Berry
bundle. Let us recall that this moduli space consists of pairs
$(S, \theta)$ where $S$ is a Bohr - Sommerfeld cycle in $M$ of the
fixed topological type and $\theta$ is a half weight on $S$ such
that the corresponding volume
$$
\int_S \theta^2 = 1
$$
--- is normalized.

For every smooth function $f \in C^{\infty}(M \to \R)$ there is
an induced smooth function $F_f$ on the moduli space which is
defined by the following formula
$$
F_f(S, \theta) = \tau \int_S f|_S \theta^2 \tag 2
$$
where $\tau$ is a real parameter. This formula (2) gives us a map
$$
C^{\infty}(M \to \R) \to C^{\infty}(\sB^{hw,1}_S \to \R). \tag 3
$$

On the other hand since the moduli space is defined in terms of
symplectic geometry one gets that every symplectomorphism $\phi$
of $(M, \om)$ induces an automorphism $\tilde{\phi}$ of
$\sB^{hw,1}_S$ which preserves all the structures $\Om, I, G$
(see [2], [3]). This gives us a dynamical correspondence
$$
\Theta_{BS}: C^{\infty}(M \to \R) \to Vect (\sB^{hw}_1)
$$
which satisfies the correspondence principle
$$
\Theta_{BS}(\{f, g \}_{\om}) = [\Theta_{BS}(f), \Theta_{BS}(g)]
$$
just due to the construction. Moreover the image $Im \Theta_{BS}$
lies in the subspace
$$
Vect_K(\sB^{hw,1}_S) \subset Vect (\sB^{hw,1}_S)
$$
consists of the vector fields which  preserve the K\"ahler
structure on the moduli space.

Now a natural question arises: for every smooth function $f$
there are two vector fields on the moduli space $\sB^{hw,1}_S$.
The first vector field is the image $\Theta_{BS}(f)$ coming from
the dynamical properties. The second is the hamiltonian vector
field $X_{F_f}$ of the induced function $F_f$. The key point is
contained in the following statement.

\proclaim{Proposition 1} For every smooth function $f$ on $M$
these two vector fields are proportional
$$
X_{F_f} = 2 \tau \Theta_{BS}(f).
$$
\endproclaim

We prove this coincidence in the next section by direct
computations. Now we would like to remark that the Main Theorem
stated in the Introduction is a consequence of this proposition.
Indeed:

1)First of all if $f$ and $g$ are two distinct functions on $M$
then the hamiltonian vector fields $X_{F_f}$ and $X_{F_g}$ are
the same iff $f = g + c$ where $c$ is a constant. But it's clear
from (2) that the constants over $M$  go to constants over
$\sB^{hw,1}_S$. Thus the map (3) is an inclusion.

2)Further, we know that $\Theta_{BS}(f)$ generates the flow which
preserves both the symplectic form $\Om$ and the riemannian metric
$G$. It means that $X_{F_f}$ preserves $G$ to.

3) And due to the fact that the dynamical correspondence
satisfies the correspondence principle one gets that
$$
\{ F_f, F_g \}_{\Om} = 2 \tau F_{\{f, g \}_{\om}}. \tag 4
$$
We have to mention here that formula (4) has been proved in [10],
[11] directly. However now the dynamical correspondence gives us
the possibility to reprove this known fact and prove new facts
listed above. As well one can see that real parameter $\tau$
plays the role of a multiple of the Planck constant.

Therefore one can see that the statement of the Main Theorem
comes directly from the dynamical correspondence. It reflects
just geometrical naturalness of the relationship between
infinitesimal symplectic deformations induced by functions $f$
and $F_f$ on symplectic manifold $M$ and K\"ahler manifold
$\sB^{hw,1}_S$ respectively.

At the same time one can see that we have just skipped the
discussion for the case of higher level $k$ (see [2]). It's easy
to see that the construction can be extended to the case of any
level. Here one wants to relate our real parameter $\tau$
appeared in the definition of $F_f$ and consequently in (4) with
$k$. Since $k$ is the inverse of the Planck constant (see f.e.
[2]) then one could take
$$
\tau = \frac{1}{2 k} $$ getting right correspondence during $BPU$
- map.

 \head \S 4. Computations
\endhead

First of all let us recall that the tangent space $T_{(S_0,
\theta_0)} \sB^{hw,1}_S$ is modeled by pairs
$$(\psi_1, \psi_2), \quad \quad \psi_i \in C^{\infty}(S_0 \to \R)
$$
such that
$$
\int_{S_0} \psi_i \theta_0^2 = 0
$$
(see [2], [3]). Now we compute the components of $\Theta_{BS}(f)$
as follows. The hamiltonian vector field $X_f$ over $M$ near
$S_0$ which is a lagrangian cycle can be decomposed into two parts
$$
X_f = V_f + W_f
$$
where
$$
V_f = \om^{-1}(d (f|_{S_0}))
$$
and
$$
W_f = X_f - V_f
$$
is parallel to $S_0$ (more rigorously $W_f$ belongs to $TS_0$).
As well one can see that $$ W_f = (\om^{-1}(df))|_{S_0}. $$
 We
understand $V_f$ as "outer" part of $X_f$ with respect to $S_0$
while $W_f$ is "inner" part whose flow preserves $S_0$. Then
$V_f$ gives us the corresponding deformation of $S_0$ itself
while $W_f$ deforms $\theta_0$ over $S_0$. It means that
$\Theta_{BS}(f)$ has the following components
$$
\aligned \psi_1(S_0, \theta_0) = f|_{S_0} - c, \\
\psi_2 (S_0, \theta_0) = \frac{Lie_{W_f} \theta_0}{\theta_0}\\
\endaligned
$$
where $c$ is the normalized constant
$$
c = \int_{S_0} f|_{S_0} \theta_0^2.
$$
The Lie derivative contained in the second equality can be
understood as follows. The square $\theta_0^2$ gives a volume
form $d\mu_0$ on $S_0$ and one can take the Lie derivative
$$
L_f = \frac{ Lie_{W_f} d \mu_0}{d \mu_0}
$$
divided by $d \mu_0$. Then
$$
Lie_{W_f}(\theta_0^2) = 2 Lie_{W_f} \theta_0 \cdot \theta_0
$$
and consequently
$$
L_f = 2 \frac{Lie_{W_f} \theta_0}{\theta_0} = 2 \psi_2 (S_0,
\theta_0).
$$
Now let us study the hamiltonian vector field $X_{F_f}$ over
$\sB^{hw,1}_S$. The differential $dF_f$ has the form
$$
dF_f(S_0, \theta_0)(\al, \be) = \tau \int_{S_0} f|_{S_0} 2 \be
\theta_0^2 + \tau \int_{S_0} d\al((\om^{-1}(df))|_{S_0})
\theta_0^2. \tag 5
$$
Due to simplicity of the symplectic form (see [2], [3]):
$$
\Om(S_0, \theta_0)<(\al, \be), (\ga, \de)> = \int_{S_0} (\al \de
- \be \ga) \theta_0^2 \tag 6
$$
one can immediately "convert" the first part
$$
\int_{S_0} f|_{S_0} \be \theta_0^2 = \Om<( \psi_1', 0), (0, \be)>
$$
such that
$$
\psi_1' = 2 \tau f|_{S_0} - c'
$$
where $c'$ is the normalized summand as above. Further, the
second summand in (5) can be rearrange as follows. First of all
$(\om^{-1}(d f))|_{S_0}$ is exactly the "inner" part $W_f$ of the
hamiltonian vector field $X_f$. Then one has
$$
\aligned \int_{S_0} d \al (W_f) \theta_0^2 = \int_{S_0} d \al
\wedge \imath_{W_f} d \mu_0 = \\ - \int_{S_0} \al d(\imath_{W_f} d
\mu_0) = - \int_{S_0} \al \frac{d \imath_{W_f} d \mu_0}{d \mu_0}
\theta_0^2\\ \endaligned
$$
where we use the integration by parts. Substituting to (6) we get
$$
\psi_2' = \tau \frac{Lie_{W_f} d \mu_0}{d \mu_0} = 2 \tau L_f.
$$
Comparing $(\psi_1, \psi_2)$ and $(\psi_1', \psi_2')$ one gets
$$
X_{F_f} = 2 \tau \Theta_{BS}(f)
$$
and the proof of the proposition is completed.

\head \S 5. ALAG - quantization
\endhead

Thus we have seen that the moduli space $\sB^{hw,1}_S$ can be
regarded as the quantum phase space for the ALAG - quantization of
given symplectic manifold $(M, \om)$. Over this infinite
dimensional K\"ahler (moreover it's an algebraic) manifold endowed
with symplectic form $\Om$, integrable complex structure $I$ and
the corresponding riemannian metric $G$ one has the space $\Cal
N$ of quantum observables consists of all smooth real functions
whose hamiltonian vector fields preserve both the symplectic
structure and the riemannian metric. Inside the space $\Cal N$
one has subspace of quantized observables
$$
\{ F_f \} \subset \Cal N
$$
which is isomorphic as a Lie algebra (up to scaling depending on
$\tau$) with $C^{\infty}(M \to \R)$. Thus with the Main Theorem
in hands we can say that we've performed kinematics and dynamical
data of the quantum theory. It remains to define

--- probabilistic (or measurement) aspects of the theory.

This question includes as well state reduction procedure. And the
problem arises in this way is based on the following fact: our
infinite dimensional K\"ahler quantum phase space is noncompact.
It implies  that some quantum observables have not critical points
(= eigenstates) at all in $\sB^{hw,1}_S$! This makes any
measurement process impossible for such observables. So one has
to construct an appropriate compactification of the moduli space
such that for any  special function $F_f$ it were a set of
critical points. There are a number of usual ways to compactify
the space (f.e. the most common is the Gel'fand approach) but
here there is some proper way which begins with the
considerations of critical points of special functions. Below we
outline this way leaving all details to [12]. At the end we
discuss two reductions of the method in cases when either real or
complex polarization is fixed over $M$.

From  formula (5) we have the following simple geometrical fact.
 \proclaim{Proposition 2} A point $(S_0, \theta_0) \in \sB^{hw,1}_S$
 is a critical one for special function $F_f$ if and only if
 the hamiltonian vector field
$X_f$ over $M$ of the original function $f$ preserves this pair:
$$
\aligned f|_{S_0} = const,\\
Lie_{W_f} \theta_0 = 0.\\
\endaligned
$$
\endproclaim

Moreover for critical set of any special function $F_f$ one has
very strong geometrical properties. If $Crit_i(F_f)$ is a
connected smooth component  of the critical set of $F_f$ then
\proclaim{Claim} The submanifold $Crit_i (F_f)$ is a complex
submanifold of $\sB^{hw,1}_S$.
\endproclaim
The proof (contained in [12]) is based on the following
description of a tangent vector to $Crit_i(F_f)$. A pair
$(\psi_1, \psi_2)$ represents the tangent vector iff
$$
Lie_{W_f} \psi_i = 0.
$$
The symmetry in the condition  means that the tangent space to
$Crit_i(F_f)$ in a smooth point is a complex subspace of $T
\sB^{hw,1}_S$.

Now the description of critical points of $F_f$ in terms
intrinsic to the symplectic geometry of the original manifold $M$
hints how one can complete the moduli space. Roughly speaking
together with lagrangian Bohr - Sommerfeld cycles one can consider
all isotropic with respect to $\om$ subcycles of fixed
topological types. Moreover one should require for these
subcycles to be intersections of lagrangian Bohr - Sommerfeld
cycles with smooth submanifolds of topological types $D, D^2, ...,
D^n$ where $D \in H_{2n-2}(M, \Z)$ is Poincare dual to  the fixed
"symplectic" class $[\om]$. Let us take a "divisor" $Y \subset M$
representing $D$ such that $Y$ is a symplectic submanifold (one
has to impose this condition as we see below). Now let
$\sB^{hw,1}_{S, Y}$ consists of all intersections $S \cap Y$ with
some corresponding half weights (one can realize this moduli
space in quite natural way: it is the moduli space of half
weighted Bohr - Sommerfeld cycles of the fixed volume over
symplectic manifold $Y$ endowed with  induced prequantization
equipment). Then there exists a smooth function $f_Y$ over $M$
such that

1) the corresponding special function $F_{f_Y}$ has not critical
points in $\sB^{hw,1}_S$ at all

but

2) has critical points in the compactification component
$\sB^{hw,1}_{S, Y}$.

It means that the induced hamiltonian flow contracts the moduli
space $\sB^{hw,1}_S$ to the boundary component $\sB^{hw,1}_{S,
Y}$. But the last one is not compact itself hence one has to
continue the process iterating to the top component which is a
"Hilbert scheme" of $M$.

Really let us take any  smooth section $s \in \Ga (L)$ whose zero
set coincides with the "divisor":
$$
(s)_0 = Y.
$$
Then
$$
f_Y = <s, s>
$$
is a nonnegative smooth function on $M$. This function is not
constant being restricted on every lagrangian cycle and this gives
us  statement 1) above. And it's not hard to check that every
point of $\sB^{hw,1}_{S, Y}$ is stabilized by the induced flow
$X_{f_Y}$. Really our function $f_Y$ is identically zero on $Y$
and moreover $Y$ is a component of the critical set $Crit(f_Y)$
of the function. It means that $d f_Y$ vanishes at each point of
$Y$ hence the same is true for the hamiltonian vector field. Thus
the corresponding action on each point $(S', \theta') \in
\sB^{hw,1}_{S, Y}$ is trivial and this remark gives us statement
2) above.

Now one can repeat the arguments for the boundary component
$\sB^{hw,1}_{S, Y}$ using a function which doesn't admit any
invariant $n-1$ - dimensional isotropic submanifold, coming to
isotropic $n-2$ - submanifolds etc.etc. and iterating the process
one gets a tower of moduli spaces compactifying the original one
where the "top" component just corresponds to sets of points of
$M$.

A compactification component $\sB^{hw,1}_{S, Y} \subset \partial
\sB^{hw,1}_S$ admits natural K\"ahler structure which is
compatible with the original structure on $\sB^{hw,1}_S$. The
same is true for each compactification component. Thus one can
expect that the compactified moduli space
$\overline{\sB}^{hw,1}_S$ is an infinite dimensional K\"ahler
manifold. The definition of special function $F_f$ can be easily
extended to every compactification component and hence to the
compactified moduli space. The dynamical correspondence again
ensures that for these extended functions all the statements of
the Main Theorem still hold. This means that one gets a real way
to quantize classical mechanical systems in the framework of ALAG
- programme. And it's more then reasonable to call this approach
as "ALAG - quantization".

\head \S 6. Polarizations
\endhead

Here we discuss two examples when ALAG - quantization can be
reduced to known ones.

\subheading{Real polarization} Let us suppose that $M, \om$ is a
completely integrable system. It means that there are exist $n$
smooth non constant function $f_i, i=1, ..., n$ in involution
$$
\{f_i, f_j \}_{\om} = 0 \quad \quad \forall i, j
$$
defining  lagrangian fibration
$$
\pi: M \to \De
$$
where $\De$ is a convex polytope in $\R^n$ and for any inner
point $t \in \De \backslash \partial \De$ the corresponding fiber
$\pi^{-1}(t)$ is a smooth lagrangian cycle. The boundary
$\partial \De$ represents degenerations: over a point of a
hyperplane in $\partial \De$ one has isotropic smooth submanifold
of dimension $n-1$ etc. (and over vortices of $\De$ one has just
points). The known quantization procedure has as the Hilbert
space $\Cal H$ in this situation the following direct sum. Namely
let us take all the lagrangian fibers which are Bohr - Sommerfeld
with respect to a prequantization quadruple over $M$. One expects
([7], [8]) that there are finitely many such fibers. Denote these
fibers as $S_i, i = 1, ..., l$. Then
$$
\Cal H = \sum_{i=1}^l \C <S_i>
$$
(for details see [7]).

What we get in this situation applying ALAG - programme? Let us
construct $\sB^{hw,1}_S$ where $S$ is the class of the fiber in
lagrangian fibration. Let us take the induced special functions
$$
F_1, ..., F_n \quad | F_i = F_{f_i}.
$$
Let us take the following intersection
$$
 \Cal P = Crit(F_1) \cap ... \cap Crit (F_n) \subset \sB^{hw,1}_S.
 $$
 Then this set is the set of mutual eigenstates of quantum
 observables $F_1, ..., F_n$. And we have the following
\proclaim{Proposition 3} The set $\Cal P$ is a double covering of
the set of Bohr - Sommerfeld fibers $\{S_i \}$.
\endproclaim

This means that one can reconstruct $\Cal H$ in terms of ALAG.
Namely, the set $\Cal P$ is a set of points. If we take the
following direct sum
$$
\Cal H_{\R} = \sum_i \R p_i
$$
then the natural antiholomorphic involution on $\sB^{hw,1}_S$
induces an involution on $\Cal P$ and consequently a complex
structure on $\Cal H_{\R}$. This complex space is canonically
isomorphic to $\Cal H$.

\subheading{Proof} We have to show that if $S_i$ is a Bohr -
Sommerfeld fiber then there exists such half weight $\theta_i$
that  pair $(S_i, \pm \theta_i)$ is a critical point of $F_j$ and
vice versa. Turning to the Proposition 2  one sees that this fact
is quite obvious --- for every $j$ our original function $f_j$ is
constant along the fiber $S_i$ so it remains to find an
appropriate invariant with respect to $W_{f_j}$ half weight
$\theta_i$. To do this firstly let us find an invariant volume
form $d \mu_i$ over $S_i$. Over $S_i$ we have $n$ nonvanishing
pointwise independent hamiltonian vector fields which are
parallel to $S_i$:
$$
X_{f_1}, ..., X_{f_n}.
$$
Let us take the corresponding set of  the differentials
$$
df_1, ..., df_n
$$
and perform the top wedge product
$$
\eta = df_1 \wedge .... \wedge df_n
$$
choosing an appropriate order. This $n$- form is totally zero
being restricted to $S_i$ but one can take a $n$ - form $\eta'$
defined by
$$
d \mu = \eta \wedge \eta' \tag 7
$$
at each point of $S_i$ where $d \mu$ is the usual symplectic
volume form. Of course this form $\eta'$ is not uniquely defined
by (7) but its restriction to $S_i$ is unique indeed. One takes
$$
d \mu_i = c \eta'|_{S_i},
$$
where $c$ is the normalizing constant. Since all forms which were
used in the description are invariant under the hamiltonian flow
induced by each $f_i$ one gets that the volume form $d \mu_i$
looks like quite canonical candidate to be the square of the half
weight which we want to construct. It remains to mention that
there are exactly two half weights over $S_i$ which give the same
form $d \mu_i$. Thus we get a pair of conjugated points $(S_i, \pm
\theta_i) \in \sB^{hw,1}_S$ which are critical for every $F_j$.
Moreover there are no other choices of invariant half weights
over $S_i$. One can check that if there exists a half weight
$\theta_i'$ which is invariant under each $f_i$ then the ratio
$$
\psi = \frac{\theta_i'}{\theta_i} \in C^{\infty}(S_i \to \R)
$$
is a smooth function which should satisfy
$$
L_{W_{f_j}} \psi = 0 \tag 8
$$
for every $j = 1, ..., n$. But again the set of hamiltonian
vector fields $W_{f_j}$ forms a complete system such that the
condition (8) implies that $\psi$ is a constant. But the
normalization condition
$$
\int_{S_i} \theta_i' = 1
$$
means that this constant is equal to $\pm 1$.

And in the other direction: if $(S_0, \theta_0)$ is a mutual
critical point for every $F_j$ then removing the half weight part
one gets a Bohr - Sommerfeld cycle $S_0$ such that for every $j$
the function $f_j$ is constant on $S_0$. But it means that $S_0$
is a fiber and the proof is completed.

\subheading{Remark} One can easily see now why the set of Bohr -
Sommerfeld fibers has to be at least discrete. Indeed, if one has
a non isolated Bohr - Sommerfeld fiber $S_0$ corresponds to an
inner point $t_0$ of $\De$
$$
t_0 \in \De \backslash \partial \De
$$
then it would be a non isolated critical point for any special
function $F_j$. This means that there exists a smooth function
$\psi \in C^{\infty}(S_0 \to \R)$ such that for any original $f_i$
one has
$$
Lie_{W_{f_i}} \psi = 0 \quad \quad i = 1, ..., n.
$$
Again using the fact that vector fields $W_{f_i}$ form complete
system over $S_0$ we get that $\psi$ has to be constant. These
arguments can be refined to the case of the faces of $\De$. The
point is that a sequence of smooth Bohr - Sommerfeld fibers can
converge a priori to a boundary point representing an isotropic
submanifold $S'$ of dimension less then $n$. Let us suppose that
this limit cycle belongs to a face of $\partial \De$ of dimension
$n-1$. Then  it has to be a point of our compactification (see
Section 5) thus one can deduce that it could not be a limit point
for such a sequence. Indeed otherwise there exists  a function
$\psi'$ over the submanifold $S'$ satisfies
$$
Lie_{W_{f_i}'} \psi' = 0, \quad \quad i = 1, ..., n-1
$$
where we change the order for $f_i$ such that the face containing
$S'$ corresponds to say maximal value of $f_n$. This function
gives us the partial deformation of $S'$ to the $n$ - dimensional
"resolution". But again this $\psi$ has to be  constant. Therefore
one has only two possibility for any sequence of regular Bohr -
Sommerfeld fibers: it is either finite or converge to a vortex of
$\partial \De$. The last case can be removed in some particular
cases.

Thus one can understand the quantization with real polarization as
a reduction of ALAG - quantization.

\subheading{Complex polarization} The usual way to quantize a
classical system using a complex polarization has been discussed
in Section 2. Here we use the approach proposed in [6] combining
together the methods of Souriau - Kostant and Berezin (see
section 2 above). So let $M, \om$ is endowed with a compatible
integrable complex structure  $I$ and we take the space of
holomorphic section $H^0(M_I, L)$ as $\Cal H$. It can be
translated (see section 1, 2  above) to the language of
projectivization. Then one has
$$\Cal P = \proj H^0(M_I, L)$$ as the phase space of quantized
system.

Again one can relate two phase spaces using so called BPU - map
(see [2], [3]):
$$
BPU: \sB^{hw,1}_S \to \Cal P.
$$
Thus the first quantum space is fibered over the second one. Now
if $f$ is a quantizable function (see [6]) then one could compare
induced quantum observables $F_f$ and $\tilde{Q}_f$ over the
quantum phase spaces. Postponing any complete computations we
just claim that the following fact takes place.
\proclaim{Proposition 4} The critical set $Crit(F_f)$ is embedded
by  $BPU$ - map to the critical set $Crit (\tilde{Q}_f$:
$$
BPU(Crit(F_f)) \subset Crit(\tilde{Q}_f).
$$
\endproclaim

The proof is quite familiar  --- we again apply the dynamical
correspondences for $\sB^{hw,1}_S$ and $\Cal P$ respectively just
noting that $BPU$ - map is invariant under any infinitesimal
deformation generating by $f$ if $f$ is a quantizable function
(recall that it means that $X_f$ preserves the complex structure
and hence acts as an infinitesimal automorphism of $\Cal P$).
Therefore due to the dynamical correspondence one has that if $f$
is a quantizable function then $$ d BPU( X_{F_f}) =
K(\tilde{Q}_f) \tag 9
$$
where $K$ denotes the Killing vector field.  This gives us the
statement together with a number of consequent remarks. We'd like
just mention here that one which makes the reduction mentioned
above. Namely if $F_f$ has an eigenstate at $(S_0, \theta_0)$ then
it implies that $\tilde{Q}_f$ has an eigenstate at
$$
p_0 = \pi (S_0, \theta_0) \in \Cal P.
$$
So again one can reduce ALAG - quantization to the known setup.
\subheading{Remark} For a generic function $f$ two critical sets
$Crit(F_f)$ and $Crit(\tilde{Q}_f$ don't coincide modulo $BPU$ -
map. The point is that in this general situation the differential
$dBPU$ "kills" the hamiltonian vector field $X_{F_f}$ in some
points covering critical points "down stair" (see formula (9)).
But if the coincidence doesn't appear for a smooth function $f$
then one can deform this function such to that for which the
coincidence $$ BPU(Crit F_f) = Crit(\tilde{Q}_f)$$ takes place.

 On
the other hand it's clear that $BPU$ - map as it is considered
now is a linearization of non linear quantum mechanical system
which is our ALAG - quantum mechanical system. At the same time
the case with real polarization above can be regarded as a
linearization as well. Thus these types of well known
quantizations are just two genuinely different linearizations of
the ALAG - problem.

Moreover one could try to exploit ALAG - quantization as a link
between two different cases of polarized mechanical system. Namely
if one has over $M, \om$ simultaneously real and complex
polarizations then $\sB^{hw,1}_S$ is an universal object which
endows us with a relationship. The constructions of the present
section give quite definite way to compare two known
quantizations. Namely if over a symplectic manifold $(M, \om)$
one has both real and complex polarizations then one can relate
two Hilbert spaces using our constructions. For this it's
sufficient to find an appropriate smooth function $f_0$ over the
original symplectic manifold such that:

1) $f_0$ is an algebraic combination of the integrals $f_1, ...,
f_n$ of the system;

2) the induced special function $F_{f_0}$ on the moduli space
$\Cal B^{hw,1}_S$ has as the critical set  precisely the set of
the points which correspond to Bohr - Sommerfield fibers;  3) the
corresponding T\"oplitz operator $T(f_0)$ on $H^0(M_I, L)$ has
pairwise different eigenvalues (and of course let $f_0$ preserves
the fixed complex structure).

Then applying Propositions 3 and 4 one gets an isomorphism between
the corresponding Hilbert spaces. Let us remark that the ambiguity
$\pm$ in taking of the half weight part is killed by $BPU$ - map.
At the same time since the critical values of $\tilde{Q}_{f_0}$
are different then the critical points of $\tilde{Q}_{f_0}$ are
isolated. It means that one has a preferred basis in $H^0(M_I,
L)$, defined by $f_0$. Proposition 4 then gives us that the images
of $(S_i, \pm \theta_i)$ belong to this critical set. Now one can
formulate the special conditions for this $f_0$ to investigate
the isomorphism. If the space of quantizable function is big
enough one can deform a given $f_0$, satisfying the conditions
1), 2), 3) above, such that Proposition 4 could be exploited in
the opposite direction thus $$ BPU(Crit F_{f}) = Crit
(\tilde{Q}_f)$$ and this would give the isomorphism.

\head Conclusion
\endhead

As it was pointed out at the end of Section 2  the qualitative
difference between classical and quantum dynamics is in presence
of an appropriate riemannian metric in the last ones
considerations. It plays crucial role in the measurement process
which is a distinguished part of the quantum theory. At the same
time for any quantum observables considered as a smooth function
one has a precise value at each quantum state which is a point of
our projective space. This fact makes some confusion in the
application of the geometrical formulation of quantum mechanics.
Thus it would be quite reasonable to derive some partial version
of what was proposed where some special "things" will play the
role of quantum observables. It would be a kind of functions on a
super manifold --- such functions have not "real" values at
points. Below we show that in the framework of ALAG - programme
one has some reduction of the ALAG - quantization to a subject
including these "things".

Our moduli space $\sB^{hw,1}_S$ is fibered over a real smooth
manifold $\sB_S$ which is the moduli space of original
"unweighted" Bohr - Sommerfeld cycles (see [2], [3]):
$$
q: \sB^{hw,1}_S \to \sB_S, \quad \quad q(S, \theta) = S.
$$
This fibration is a Lagrangian fibration. Every smooth function
$f \in C^{\infty}(M \to \R)$ defines a vector field $Y_f$ on
$\sB_S$ due to the following arguments. At each $S \in \sB_S$ the
restriction
$$
\psi = f|_S
$$
gives us a tangent vector
$$
\psi \in T_S \sB_S = C^{\infty}(S \to \R)/const.
$$
Generalizing this remark over whole $\sB_S$ one gets the vector
field $Y_f$. On the other hand as usual we have a dynamical
correspondence
$$
\Theta_S: C^{\infty}(M \to \R) \to Vect (\sB_S)
$$
due to the fact that every symplectomorphism of $M, \om$ acts as
an automorphism of $\sB_S$. It's a weaker version of Proposition 2
which states that
$$
Y_f = \Theta_S(f).
$$
It implies that map $\Theta_S$ preserves the Lie algebra
structure.

For every function $f$ one has not any value of this function on
a Bohr - Sommerfeld cycle $S \in \sB_S$ unless $f$ is constant
alone the submanifold. It means that one has a "measurement"
process for $\tilde{Y}_f$ over $\sB_S$ in points where the vector
field $Y_f$ vanishes. Here we denote as $\tilde{Y}_f$ the "thing"
which we would like to define. For these "things" we have a
Poisson bracket:
$$
\{ \tilde{Y}_f, \tilde{Y}_g \}_B = \widetilde{[Y_f, Y_g]}.
$$
It's quite easy to check that this one satisfies all usual
conditions. It can be done using the following obvious equality
$$
\{ \tilde{Y}_f, \tilde{Y}_g \}_B = \widetilde{Y_{\{f, g
\}}}_{\om}.
$$

Now we  can recognize what are these "things" which correspond to
smooth real function over $M$. Namely (see f.e. [4]) each
$\tilde{Y}_f$ looks like a function on an odd supersymplectic
manifold. This manifold is constructed over $\sB$ using usual
procedure (see [4]) which "twist" the standard symplectic
structure on $T^* \sB_S$ getting an example of the
supersymplectic structure on
$$
 \Pi T^* \sB_S.
 $$
  Functions in this setup are represented by multivector fields
  --- as well as $\tilde{Y}_f$ has the representation by
  $Y_f$ over $\sB_S$ for every smooth function $f$. Then one can
  see that the Poisson bracket $\{.,.\}_B$ defined above is
  nothing but so called Buttin bracket (see [4]). Our "quantized
  observables" $\tilde{Y}_f$ are distinguished by the condition
  that they are represented by vector fields (so they are of
  degree 1 in odd variables).

  We are not ready to go further, leaving this story with super
  symplectic geometry. We just have to mention that ALAG -
  programme is in some sense based on super symplectic geometry.
  Really one could derive that the notion of Bohr - Sommerfeld
  lagrangian cycle is absolutely equivalent to the notion of
  lagrangian cycle in some appropriate even super symplectic
  manifold. This super symplectic manifold is nothing but the
  principle bundle $P \to M$ associated to our prequantization
  hermitian line bundle $L \to M$ endowed with our prequantization
  connection $a$. Over this principle bundle $P$ one has
  an even super symplectic form $\om_e$ defined as follows:
  for every pair of tangent to $P$ vectors $v_1, v_2$ in a point
  $p \in P$ one can decompose these ones with respect to the
  fixed connection into horizontal and vertical part. Then
  $\om_e$ acts on the horizontal parts as usual $\om$ (in skew
  symmetric style) while on
the vertical parts as the riemannian pairing (in symmetric style).
Now one can say what is a Planckian cycle $\tilde{S} \in \Cal P$
(the definition see in [2]): it is just a lagrangian cycle with
respect to even super symplectic form $\om_e$.

Thus one could say that ALAG - programme is based on the super
symplectic geometry. It means that one should expect some new
results comparing ALAG - constructions with constructions belong
to this super symplectic geometry as well as a generalization of
ALAG ("non abelian lagrangian algebraic geometry") in the setup
of super symplectic geometry.

\enddocument